\documentclass{amsart}
\usepackage{graphicx}
\theoremstyle{definition}

\theoremstyle{remark}

\numberwithin{equation}{section}

\begin{document}

\title[]{On the equivariant algebraic Jacobian for curves of genus two}
\author{Chris Athorne}%
\address{School of Mathematics and Statistics\\University of Glasgow\\Glasgow G12 8QW\\UK}%
\email{christopher.athorne@gla.ac.uk}%

\thanks{}%
\subjclass{}%
\keywords{}%

\begin{abstract}
We present a treatment of the algebraic description of the Jacobian
of a generic genus two plane curve which exploits an $SL_2(k)$
equivariance and clarifies the structure of E.V.Flynn's 72 defining
quadratic relations. The treatment is also applied to the Kummer
variety.
\end{abstract}
\maketitle
\section{Introduction}

The work described in this paper is a reflection on some material in Chapters 2
and 3 of \cite{CF}. We intend to present a simplification of
the explicit description of the algebraic Jacobian , ${\mathcal
J(C)}$, for a genus 2 curve given there and at \cite{FlynnWeb}.

The Jacobian of a non-singular, compact Riemann surface, $\mathcal
X$, is the group $Pic^0$ of divisors of degree zero factored out by
principle divisors. This can be constructed analytically using the
Abel map \cite{FK}. As such, $g$ being the genus of $\mathcal C$,
the Jacobian is $\mathbb C^g/\Lambda$, $\Lambda$ being the
$g$-dimensional lattice of periods.

The Riemann-Roch theorem describes the dimensions of linear spaces
of functions with prescribed poles on $\mathcal X$. These functions
are coordinates on $\mathcal X$ and relations between them provide
us with (generally singular) models of the surface as an algebraic
curve in some projective space. An economical description is
obtained by taking $P\in\mathcal X$ to be a Weierstra{\ss} point. In
the case of the genus 2 (hyperelliptic) surface, there are
coordinates $x:{\mathcal X}\mapsto \mathbb P^1$ and $y:{\mathcal
X}\mapsto \mathbb P^1$ with poles of orders 2 and 5 respectively which
satisfy a relation of the form
\[y^2=4x^5+\lambda_4x^4+\lambda_3x^3+\lambda_2x^2+\lambda_1x+\lambda_0,\]
the $\lambda_i$ being constants in the ground field. As a curve in
$\mathbb P^2$ this is singular at infinity. Functions associated
with more general special divisors provide us with other models.
Such models are related by birational transformations. Thus we will
be concerned, as is \cite{CF}, with (singular) models of the genus 2
curve in the form
\[y^2=g_6x^6+6g_5x^5+15g_4x^4+20g_3x^3+15g_2x^2+6g_1x+g_0,\]
which are related amongst themselves and to the quintic by simple
M\"obius maps:
\[x\mapsto\frac{\alpha x+\beta}{\gamma x+ \delta},\]
\[y\mapsto\frac{y}{(\gamma x+ \delta)^3}.\]

A similar philosophy allows the algebraic construction of ${\mathcal
J(C)}$ and such constructions are of use over fields other than
$\mathbb C$ which is what partly motivates \cite{CF}.

Another construction goes back to Jacobi and is described in
\cite{Mumtata}.

In \cite{CF} $Pic^0$ is identified with $Pic^2$ and ${\mathcal
J(C)}$ is constructed as a quadric variety in $\mathbb P^{15}$, the
locus of seventy two linearly independent quadratic identities.
Sixteen homogeneous coordinates on $\mathbb P^{15}$ are chosen to be
symmetric functions in two points on the curve: in our notation,
$(x_1,y_1)$ and $(x_2,y_2)$. These coordinates are allowed to have
poles of order up to 4 on the special divisor $D=(x,y)+(x,-y),$ that
is, when $x_1=x_2$ and $y_1=-y_2$ and up to order 2 at the
(singular) point at infinity. By looking at quadratic expressions in
these coordinates and by balancing poles \cite{CF} construct the
seventy two identities to be found explicitly at \cite{FlynnWeb}

The purpose of the current paper is to use a little representation
theory to oil the wheels of this machinery and to uncover some
structure intrinsic to the collection of quadratic identities. Such
an approach has already proven valuable in the analytic context
\cite{Ath1,Ath2} following upon the work of \cite{BEL}.

The idea is that the coordinates on $\mathcal J(C)$ can be chosen to
belong to irreducible $G$-modules where $G$ is a group of birational
transformations. Quadratic functions arise by tensoring up these
modules and decomposing into irreducibles. It suffices to work only
with highest weight elements and it turns out that the dimensions of
the components of the decomposition are graded by degrees of poles
on the divisor in a ``helpful" way. Because of this we need do less
work and the identities are arranged for us by the representation
theory into patterns.

We use the group $PSL_2(k)$ defined by the M\"obius transformations
above and take $k$ to be an algebraically closed field of zero
characteristic so that we can stay close to classical representation
theory as presented in, say, \cite{H}

In fact \cite{CF} does mention that their coordinates transform
nicely under translations and inversion of $x$ and goes so far as to
write down a more complicated basis which would presumably be
similar to that we present below. But the authors do not pursue the observation.

In the next section we define our notation, presenting the Lie
algebraic action of the coordinate transformations on the variables
and the coefficients of the curve and we define the construction of
a highest weight element that we use for a component of the
decomposition. We define our \emph{inhomogeneous} coordinates on
$\mathbb P^{15}$ and classify them according to dimension and degree
of pole divisor. The inhomogeneous coordinates seem easier to use
with the $PSL_2(k)$ action.

We give the decompositions of tensor products of coordinates
according to dimension and divisor degree and indicate how the
strategy of balancing dimensions and poles works to create quadratic
identities. The quadratic identities themselves we summarise in the
next section. We also need to construct various invariants and
covariants out of the coefficients of the curve which themselves are
a basis for a 7 dimensional $SL_2(k)$ module.

Some of the algebraic manipulation is done by hand. Some of it is
best done using a computer algebra package such as MAPLE, used in
this instance. Since one knows exactly where to look for
cancelations the only issue is calculating the coefficients, a
matter of linear algebra.

After that we give a construction of the Kummer variety associated
to the genus 2 curve which is perhaps a little different to the
\cite{CF} construction, though the result is entirely,
equivariantly, equivalent.

Finally we make some closing remarks.

\section{Notation}
\subsection{Modules and Tensor products}
We use the following normalisation for a basis
$\{v_0,v_1,\ldots,v_{n-1}\}$ of a standard irreducible ${\mathfrak
sl}_n$-module, $V_n$ of dimension $n$:
\begin{eqnarray}
e(v_i)&=&(n-i)v_{i-1}\nonumber\\
e(v_0)&=&0\nonumber\\
f(v_i)&=&(i+1)v_{i+1}\nonumber\\
f(v_{n-1})&=&0\nonumber\\
h(v_i)&=&(n-2i-1)v_i\nonumber
\end{eqnarray}
for $i=0,\ldots n-1.$

We call $v_0\in \ker e$ a highest weight element.

The coefficients $g_i$ of the curve are a basis for a seven
dimensional dual module so, for the sake of uniformity, we introduce
a new set of coefficients \[g^*_i=(-1)^i{\binom{6}{i}}g_{6-i}\quad
i=0,\ldots 6\] carrying the standard representation.

The tensoring of modules of dimensions $n$ and $m$, $n\geq m$ leads
to a decomposition of the form
\[V_n\otimes V_m\simeq\bigoplus_{i=n-m+1}^{n+m-1}V_i\]
and we construct the highest weight elements of the components in
this plethysm according to the following rule:

\[{(U_n\otimes V_m)}_{n+m-p,0}=\sum_{i=0}^{p}(-1)^i\frac{(n-i-1)!}{(n-1)!}\frac{(m-p+1-1)!}{(m-1)!}u_iv_{p-i}\]

The basis elements of the representations we use are to be
inhomogeneous coordinate functions defined on the Jacobian of the
genus two curve in $\mathbb P^{15}.$ Each has a singularity on the
diagonal $x_1=x_2,\,y_1=y_2$ and on the divisor. All the elements of a given irreducible
have, in fact, the same singularities as can be verified by
inspection of all the terms to be defined. We denote by ${\bf
n}^{div}_{diag}$ the class of $n$ dimensional modules with poles of
order $diag$ (respectively $div$) on the diagonal (respectively
divisor) of the product $\textsl{C}\times\textsl{C}$.

\subsection{Fundamental Irreducibles}

Let $\Delta=x_1-x_2$. There is an invariant,  related to the polar
form of the curve, namely:
\[{\mathcal I}=\frac{F(x_1,x_2)-y_1y_2}{\Delta^3}\in {\bf 1}^{1}_{3}\]
where $F(x_1,x_2)$ is the (equivariant) polar form:
\begin{eqnarray}
F(x_1,x_2)&=&g_0+3(x_1+x_2)g_1+3(x_1^2+3x_1x_2+x_2^2)g_2\nonumber\\
&&+(x_1^3+9x_1^2x_2+9x_1x_2^2+x_2^3)g_3+3x_1x_2(x_1^2+3x_1x_2+x_2^2)g_4\nonumber\\
&&+3x_1^2x_2^2(x_1+x_2)g_5+x_1^3x_2^3g_6\nonumber
\end{eqnarray}

The polar form plays a fundamental role in all approaches to the
theory. In the analytic description of the Jacobian the \emph{equivariant}
polar form allows the construction of \emph{equivariant} $\wp$-functions
\cite{Ath1}.

We introduce the set of inhomogeneous projective coordinates on
$\mathbb P^{15}$ designated by lists of standard basis elements,
${\bf P}({\bf n})^p_q\in{\bf n}^p_q:$

\begin{eqnarray}
{\bf P}({\bf
5})_{2}^{2}&=&(\frac{1}{\Delta^2},\frac{2(x_1+x_2)}{\Delta^2},
\frac{x_1^2+4x_1x_2+x_2^2}{\Delta^2},\frac{2x_1x_2(x_1+x_2)}{\Delta^2},\frac{x_1^2x_2^2}{\Delta^2})\nonumber\\
{\bf P}({\bf
4})_{2}^{3}&=&(\frac{y_1-y_2}{\Delta^3},\frac{3(x_2y_1-x_1y_2)}{\Delta^3},\frac{3(x_2^2y_1-x_1^2y_2)}
{\Delta^3},\frac{x_2^3y_1-x_1^3y_2}{\Delta^3})\nonumber\\ {\bf
P}({\bf 3})^{4}_{2}&=&(2\frac{\mathcal
I}{\Delta},\frac{(x_1+x_2)\mathcal I}{\Delta},2\frac{x_1x_2\mathcal
I}{\Delta})\nonumber\\ {\bf P}({\bf
2})^{5}_{2}&=&(\frac{y_1{\mathcal I},_{x_1}+y_2{\mathcal
I},_{x_2}}{\Delta},\frac{x_2y_1{\mathcal I},_{x_1}+x_1y_2{\mathcal
I},_{x_2}}{\Delta})\nonumber\\ {\bf P}({\bf 1})^{6}_{2}&=&{\mathcal
I}^2\nonumber
\end{eqnarray}

These coordinates are inspired by the sixteen homogeneous coordinate
functions of Cassels and Flynn \cite{CF} but have one crucial
difference.  The fifteen inhomogeneous coordinates, created by dividing through by the coordinate $(x_1-x_2)^2$ in \cite{CF}, above have been
adapted to a decomposition of $\mathbb P^{15}$ into irreducible
$\mathfrak sl_2$-modules:

\[{\bf 15}\simeq{\bf 5}^2_2\oplus{\bf 4}^3_2\oplus{\bf 3}^4_2\oplus{\bf 2}^5_2\oplus{\bf 1}^6_2\]

It is important to notice that there is a pole grading on the
divisor related to the module dimension. As inhomogeneous coordinates
they are well-behaved at infinity and have poles of orders up to
$4+2=6$ on the divisor.

\section{Tensor products and pole gradings}
Quadratic functions on the Jacobian arise by tensoring up the
coordinates. In the (symmetric) table below are summarised the
irreducible decompositions of the symmetric tensor products (denoted $\odot$) of all
the coordinate modules.

\[
\begin{array}{ccccccc}
\odot &\vline& {\bf P}({\bf 5})^{2}_{2} & {\bf P}({\bf 4})^{3}_{2}
&{\bf P}({\bf 3})^{4}_{2} &{\bf P}({\bf 2})^{5}_{2} & {\bf P}({\bf 1})^{6}_{2}\\
&\vline&&&&&\\
\hline
&\vline&&&&&\\
{\bf P}({\bf 5})^{2}_{2} &\vline& {\bf 9}^{4}_{4}\oplus{\bf 5}^{2}_{2}\oplus{\bf 1}^{0}_{0} & & & & \\
&\vline&&&&&\\
{\bf P}({\bf 4})^{3}_{2} &\vline&{\bf 8}^{5}_{4}\oplus{\bf 6}^{3}_{4}\oplus{\bf 4}^{3}_{2} & {\bf 7}^{6}_{4}\oplus{\bf 3}^{4}_{4}& & & \\
&\vline&&&&&\\
{\bf P}({\bf 3})^{4}_{2} &\vline& {\bf 7}^{6}_{4}\oplus{\bf 3}^{4}_{2}&{\bf 6}^{7}_{4}\oplus{\bf 4}^{5}_{4} &{\bf 5}^{8}_{4}\oplus{\bf 1}^{6}_{2} & & \\
&\vline&&&&&\\
{\bf P}({\bf 2})^{5}_{2} &\vline&{\bf 6}^{7}_{4}\oplus{\bf 4}^{5}_{4} &{\bf 5}^{8}_{4}\oplus{\bf 3}^{6}_{4} & {\bf 4}^{9}_{4}\oplus{\bf 2}^{7}_{4}& {\bf 3}^{10}_{4}& \\
&\vline&&&&&\\
{\bf P}({\bf 1})^{6}_{2} &\vline&{\bf 5}^{8}_{4} &{\bf 4}^{9}_{4} &{\bf 3}^{10}_{4} &{\bf 2}^{11}_{4} & {\bf 1}^{12}_{4}\\
\end{array}
\]

\vspace{.2in}

We make some remarks about these decompositions.

Firstly the highest dimensional component in each decomposition has
predictable singularity structure: the highest weight element is
simply the product of the highest weight elements in each factor.
More generally the pole and dimension grading are connected in an
interesting but currently obscure manner.

Secondly there are ``holes" in the table. We expect to find, for
example, a $\bf 2$ inside ${\bf 5}\odot{\bf 4}.$ Its absence is due
to its vanishing. Such cancelations are the simplest of the
quadratic identites which will describe the Jacobian as a locus in
$\mathbb P^{15}.$ The identities arising in this way are nine in
number and shown in section (\ref{trivial}.

We will use the notation ${[}\cdot{]}_{\bf n}$ to denote projection
onto the $n$-dimensional irreducible of the decomposition.

In the next section we will summarise the seventy one quadratic
relations in such a concise form but before doing so we explain how
they are derived by considering the most complicated of them.

Consider ${\bf P}({\bf 2})^5_2\odot{\bf P}({\bf 2})^5_2\in{\bf
3}^{10}_4.$ The only possible cancelation is with ${\bf P}({\bf
3})^4_2\odot{\bf P}({\bf 1})^6_2.$ We take an arbitrary linear
combination of the highest weight elements, look at the worst
singularity on the divisor and choose the (one) free parameter to
kill it. This leaves us with a highest weight element in ${\bf
3}^8_4.$ Such singularites do not occur in the table except as ${\bf
5}^8_4$. We can form elements of ${\bf 3}^8_4$ by tensoring up the
the three occurances of ${\bf 5}^8_4$ with the seven dimensional
module of degree one in the coefficients of the curve, $g^*_i$, denoted $\bf g$.
Choosing the free parameters in a linear combination appropriately
we can cancel down to a highest weight in ${\bf 3}^6_4$. We continue
this process systematically until we reach a vanishing element and
we are done. In the current instance it becomes necessary to use
modules arising from tensor products of degrees two and three in the
curve coefficients. We summarise the structure of such
representations in the next section, giving their highest weight
elements before presenting the full list of quadratic identities.

\section{Quadratic relations}
\subsection{Irreducibles in the curve coefficients} Standard partition
counting arguments \cite{FH} yield plethysms for symmetric tensor products of the $V_n$.
From the decomposition $V_7\odot V_7\simeq V_{13}\oplus V_{9}\oplus
V_5\oplus V_1$ we obtain the following highest weight elements for
quadratic representations.
\begin{eqnarray}
{[}{\bf g}\odot {\bf g}{]}_{{\bf 13},0}&=&g_6^2\nonumber\\
{[}{\bf g}\odot {\bf g}{]}_{{\bf 9},0}&=&g_6g_4-g_5^2\nonumber\\
{[}{\bf g}\odot {\bf g}{]}_{{\bf 5},0}&=&\frac{1}{2^2.3}(g_6g_2-4g_5g_3+3g_4^2)\nonumber\\
{[}{\bf g}\odot {\bf g}{]}_{{\bf
1},0}&=&\frac{1}{2^3.3^2.5}(g_6g_0-6g_5g_1+15g_4g_2-10g_3^2)\nonumber
\end{eqnarray}
The cubic irreducibles, $V_7\odot V_7\odot V_7\simeq V_{19}\oplus
V_{15}\oplus V_{13}\oplus V_{11}\oplus V_9\oplus V_7\oplus V_7\oplus
V_3,$ have the following highest weight elements. Note there are two
seven dimensional irreducibles.
\begin{eqnarray}
{[}{\bf g}\odot {\bf g}\odot {\bf g}{]}_{{\bf 19},0}&=&g_6^3\nonumber\\
{[}{\bf g}\odot {\bf g}\odot {\bf g}{]}_{{\bf 15},0}&=&\frac{8}{11}g_6(g_6g_4-g_5^2)\nonumber\\
{[}{\bf g}\odot {\bf g}\odot {\bf g}{]}_{{\bf 11},0}&=&\frac{3}{2.11}(g_6^2g_3+3g_6g_5g_4-2g_5^3)\nonumber\\
{[}{\bf g}\odot {\bf g}\odot {\bf g}{]}_{{\bf 9},0}&=&\frac{13}{2^4.3^2.11}(g_5g_0^2-5g_4g_1g_0+2g_3g_2g_0+8g_1^2g_3-6g_1g_2^2)\nonumber\\
{[}{\bf g}\odot {\bf g}\odot {\bf g}{]}_{{\bf 7},0}&=&\frac{1}{2^5.3^3.5.7.11}(-2778g_5g_1g_0+3795g_4g_2g_0+3150g_1^2g_4\nonumber\\
&&-1480g_3^2g_0-6300g_1g_2g_3+3150g_2^3+463g_0^2g_6)\nonumber\\
{[}{\bf g}\odot {\bf g}\odot {\bf g}{]}'_{{\bf 7},0}&=&\frac{1}{2^4.3.5.7}(-6g_5g_1g_0+165g_4g_2g_0-150g_1^2g_4-160g_3^2g_0\nonumber\\
&&+300g_1g_2g_3-150g_2^3+g_0^2g_6)\nonumber\\
{[}{\bf g}\odot {\bf g}\odot {\bf g}{]}_{{\bf
3},0}&=&\frac{1}{2^3.3^2.7}(-3g_5g_3g_0+3g_5g_1g_2+2g_4^2g_0-g_4g_3g_1-3g_4g_2^2\nonumber\\
&&+2g_2g_3^2+g_2g_0g_6-g_1^2g_6)\nonumber
\end{eqnarray}
\vspace{.1in}

Cancelations of poles occur between linear combinations of quadratic
expressions at each pole order. In the following paragraphs we list
the identities by pole order and to simplify notation use $\bf n$ for ${\bf P}({\bf n})$.

\subsection{$(\cdot)^0_0$ 10 identities}\label{trivial}
\[{[}\,{\bf 5}\odot{\bf 5}\,{]}_{\bf 1}=\frac{1}{2^4.3^2}\]
\[{[}\,{\bf 5}\odot{\bf 4}\,{]}_{\bf 2}=0\]
\[{[}\,{\bf 5}\odot{\bf 3}\,{]}_{\bf 5}=0\]
\[{[}\,{\bf 4}\odot{\bf 3}\,{]}_{\bf 2}=0\]
\subsection{$(\cdot)^2_2$ 5 identities}
\[{[}\,{\bf 5}\odot{\bf 5}\,{]}_{\bf 5}+\frac{1}{2^2.3}{\bf 5}=0\]
\subsection{$(\cdot)^3_2$ 4 identities}
\[{[}\,{\bf 5}\odot{\bf 4}\,{]}_{\bf 4}-\frac{1}{2^2.3}{\bf 4}=0\]
\subsection{$(\cdot)^4_2$ 3 identities}
\[{[}\,{\bf 5}\odot{\bf 3}\,{]}_{\bf 3}+\frac{1}{2.3}{\bf 3}=0\]
\subsection{$(\cdot)^5_2$ 2 identities}
\[{\bf 2}-\left[\,2^5.3^3.5\,{\bf g}\odot[{\bf 5}\odot{\bf 4}]_{\bf 8}+\frac{2^4.3^3.5}{7}\,{\bf g}\odot[{\bf 5}\odot{\bf 4}]_{\bf 6}\,\right]_{\bf
2}=0\]
\subsection{$(\cdot)^6_2$ 2 identities}
\[{[}\,{\bf 3}\odot{\bf 3}\,{]}_{\bf 1}+\frac{1}{2^2}{\bf 1}=0\]
\[{\bf 1}-\left[\,2^4.3^2.5\,{\bf g}\odot[{\bf 4}\odot{\bf 4}]_{\bf 7}-2^7.3^4.5.7\,[{\bf g}\odot{\bf g}]_{\bf 9}\odot[{\bf 5}\odot{\bf 5}]_{\bf 9}-\frac{2^8.3^3}{7}\,[{\bf g}\odot{\bf g}]_{\bf 5}\odot{\bf 5}-108\,{\bf g}\odot{\bf g}\,\right]_{\bf
1}=0\]
\subsection{$(\cdot)^3_4$ 0 identities}
\subsection{$(\cdot)^4_4$ 3 identities}
\[\left[\,{\bf 4}\odot{\bf 4}+2^4.3^2.5\,({\bf g}\odot[{\bf 5}\odot{\bf
5}]_{\bf 9})-\frac{2^2.3^2}{7}\,({\bf g}\odot{\bf 5})\,\right]_{\bf
3}-{\bf 3}=0\]
\subsection{$(\cdot)^5_4$ 8 identities}
\[\left[\,{\bf 4}\odot{\bf
3}+2^3.3^3.5\,{\bf g}\odot[{\bf 5}\odot{\bf 4}]_{\bf
8}-\frac{2^4.3^4}{7}\,{\bf g}\odot[{\bf 5}\odot{\bf 4}]_{\bf
6}+\frac{2.3^3}{5}\,{\bf g}\odot{\bf 4}\right]_{\bf 4}=0\]
\[\left[{\bf 5}\odot{\bf 2}-2^4.3^2.5\,{\bf g}\odot[{\bf 5}\odot{\bf 4}]_{\bf 8}+
\frac{2^5.3^3}{7}\,{\bf g}\odot[{\bf 5}\odot{\bf 4}]_{\bf
6}+\frac{2^3.3}{5}\,{\bf g}\odot{\bf 4}\,\right]_{\bf 4}=0\]
\subsection{$(\cdot)^6_4$ 10 identities}
\[\left[\,{\bf 4}\odot{\bf 2}+2^3.3^2.5\,{\bf g}\odot[{\bf 4}\odot{\bf 4}]_{\bf 7}\,\right]_{\bf 3}=0\]
\[\left[\,{\bf 4}\odot{\bf 4}-2\,{\bf 5}\odot{\bf 3}-2^3.3^3\,{\bf g}\odot[{\bf 5}\odot{\bf 5}]_{\bf 9}-\frac{2^4.3}{7}\,{\bf g}\odot{\bf 5}-\frac{3}{2.5}\,{\bf g}\,\right]_{\bf 7}=0\]
\subsection{$(\cdot)^7_4$ 8 identities}
\[\left[\,{\bf 3}\odot{\bf 2}-2^2.3^2.5\,{\bf g}\odot[{\bf 4}\odot{\bf 3}]_{\bf 6}\,\right]_{\bf 2}=0\]
\[\left[\,{\bf 5}\odot{\bf 2}-3\,{\bf 4}\odot{\bf 3}+2^4.3^2\,{\bf g}\odot[{\bf 5}\odot{\bf 4}]_{\bf 8}-\frac{2^3.3^4}{7}\,{\bf g}\odot[{\bf 5}\odot{\bf 4}]_{\bf 6}+\frac{2.3}{5}\,{\bf g}\odot{\bf 4}\,\right]_{\bf 6}=0\]
\subsection{$(\cdot)^8_4$ 10 identities}
\[{[}\,{\bf 3}\odot{\bf 3}-{\bf 5}\odot{\bf 1}\,{]}_{\bf 5}=0\]
\[\left[\,{\bf 4}\odot{\bf 2}-6\,{\bf 5}\odot{\bf 1}-2^3.3^2.5\,{\bf
g}\odot[{\bf 5}\odot{\bf 3}]_{\bf 7}\right.\]
\[\left.-2^5.3^5.5\,[{\bf g}\odot{\bf g}]_{\bf 9}\odot[{\bf 5}\odot{\bf
5}]_{\bf 9} +\frac{2^5.3^6}{7}\,[{\bf g}\odot{\bf g}]_{\bf
5}\odot[{\bf 5}\odot{\bf 5}]_{\bf 9}\right.\]
\[\left.+\frac{2^4.3^3}{7}\,[{\bf g}\odot{\bf g}]_{\bf 9}\odot{\bf
5}-\frac{2^4.3^5.5}{7^2}\,[{\bf g}\odot{\bf g}]_{\bf 5}\odot{\bf
5}+2^3.3^4\,[{\bf g}\odot{\bf g}]_{\bf 1}\odot{\bf
5}+\frac{2^2.3^2}{5}\,{\bf g}\odot{\bf g}\,\right]_{\bf 5}=0\]
\subsection{$(\cdot)^9_4$ 4 identities}
\[\left[\,{\bf 3}\odot{\bf 2}-3\,{\bf 4}\odot{\bf 1}+2^4.3^2\,{\bf g}\odot[{\bf 4}\odot{\bf 3}]_{\bf 6}+
\frac{2^2.3^3}{5}\,{\bf g}\odot[{\bf 4}\odot{\bf 3}]_{\bf
4}\,\right]_{\bf 4}=0\]
\subsection{$(\cdot)^{10}_4$ 3 identities}
\[\left[\,{\bf 2}\odot{\bf 2}-2.3^2\,{\bf 3}\odot{\bf
1}-\frac{2^5.3^5}{5}\,{\bf g}\odot[{\bf 5}\odot{\bf 1}]_{\bf
5}+\frac{2^2.3^2.31}{5}\,{\bf g}\odot[{\bf 4}\odot{\bf 2}]_{\bf
5}\right.\]
\[ -2^4.3^4.5^2.7\,[{\bf g}\odot{\bf g}]_{\bf
9}\odot[{\bf 5}\odot{\bf 3}]_{\bf 7}-\frac{2^4.3^4}{5}\,[{\bf
g}\odot{\bf g}]_{\bf 5}\odot[{\bf 4}\odot{\bf 4}]_{\bf
3}-2^3.3^5\,[{\bf g}\odot{\bf g}]_{\bf 1}\odot{\bf 3}\]
\[-\frac{2^{11}.3^7.11}{5}\,[{\bf g}\odot{\bf g}\odot{\bf
g}]_{\bf 11}\odot[{\bf 5}\odot{\bf 5}]_{\bf 9}
+\frac{2^9.3^8.7.11}{5.13}\,[{\bf g}\odot{\bf g}\odot{\bf g}]_{\bf
9}\odot[{\bf 5}\odot{\bf 5}]_{\bf 9}\]
\[+\frac{2^7.3^6.11.61}{5}\,[{\bf g}\odot{\bf g}\odot{\bf g}]_{\bf
7}\odot[{\bf 5}\odot{\bf 5}]_{\bf
9}-\frac{2^6.3^6.23.211}{5.7}\,[{\bf g}\odot{\bf g}\odot{\bf
g}]'_{\bf 7}\odot[{\bf 5}\odot{\bf 5}]_{\bf 9}\]
\[+\frac{2^9.3^7.11.157}{5^2.7}\,[{\bf g}\odot{\bf g}\odot{\bf g}]_{\bf
7}\odot[{\bf 5}\odot{\bf 5}]_{\bf
5}-\frac{2^{16}.3^6.11.13}{5^2.7^2}\,[{\bf g}\odot{\bf g}\odot{\bf
g}]'_{\bf 7}\odot[{\bf 5}\odot{\bf 5}]_{\bf 5}\]
\[\left.-\frac{2^9.3^6.11}{5^2}\,[{\bf g}\odot{\bf g}\odot{\bf
g}]_{\bf 3}\odot[{\bf 5}\odot{\bf 5}]_{\bf
5}+\frac{2^5.3^6.7}{5^2}\,{\bf g}\odot{\bf g}\odot{\bf
g}\,\right]_{\bf 3}=0\]

These relations are all linearly independent. This is guaranteed by
the pole grading and by checking within each graded component. Thus
at $(\cdot)^6_4$, for example, we have possible cancelations between poles in
$[{\bf 4}\odot{\bf 4}]_{\bf 7}$, $[{\bf 5}\odot{\bf 3}]_{\bf 7}$ and
$[{\bf 4}\odot{\bf 2}]_{\bf 3}$. Seven relations obtain from
canceling $[{\bf 4}\odot{\bf 4}]_{\bf 7}$ against $[{\bf 5}\odot{\bf
3}]_{\bf 7}$. These cannot then involve $[{\bf 4}\odot{\bf 2}]_{\bf
3}$.Three relations come from expressing $[{\bf 4}\odot{\bf 2}]_{\bf
3}$ in terms of a tensor product of ${\bf g}$ with either $[{\bf
4}\odot{\bf 4}]_{\bf 7}$ or $[{\bf 5}\odot{\bf 3}]_{\bf 7}$. But the
difference of these two possibilities arises exactly from tensoring
the seven former identities with ${\bf g}$. Hence there are exactly
7+3=10 linearly independent relations.

\section{The Kummer variety}
The Kummer is a simple, quartic relation in $\mathbb P^3$ which
contains important information about the Jacobian. It is a degree
four homogeneous relation between four of the above coordinates on
$\mathbb P^{15}.$ In \cite{CF} the Kummer relation is expressed in
terms of the variable set $(1,x_1+x_2,x_1x_2,(x_1-x_2)\tilde
{\mathcal I})$ where $\tilde {\mathcal I}$ is a non-equivariant
version of $\mathcal I.$ From the equivariant point of view the
appropriate variables are $\bf 3$ and $\bf 1.$

So we seek the corresponding quartic relation between the one- and
three-dimensional irreducible parts of $\mathbb P^{15}.$ Either by
eliminating ${\bf 1}$ and $[{\bf 4}\odot{\bf 4}]_{\bf 7}$ between
relations belonging to $(\cdot)^6_2$ and $(\cdot)^6 _4$ above or by
employing the same methodology we used to obtain them, we find:
\[\left[{\bf 3}\odot{\bf 3}+2^3.3^2.5\,{\bf g}\odot[{\bf 5}\odot{\bf
3}]_{\bf 7}+2^5.3^4.5.7\,[{\bf g}\odot{\bf g}]_{\bf 9}\odot[{\bf
5}\odot{\bf 5}]_{\bf 9}\right.\]\[\left.+\frac{2^6.3^3}{7}\,[{\bf
g}\odot{\bf g}]_{\bf 5}\odot{\bf 5}+3^3\,{\bf g}\odot{\bf
g}\right]_{\bf 1}=0,\] to which we may add the relations already
found above, under $(\cdot)^8_4$,
\[{[}\,{\bf 3}\odot{\bf 3}-{\bf 5}\odot{\bf 1}\,{]}_{\bf 5}=0.\]

Since tensoring by the invariant $\bf 1$ is simply ordinary
multiplication we can eliminate all occurrences of $\bf 5$ by
multiplying the first relation by ${\bf 1}^2$ to obtain an
invariant, quartic, homogeneous expression in $\bf 3$ and $\bf
1$:

\[[{\bf 3}^2]_{\bf 1}.{\bf 1}^2+2^3.3^2.5\,{\bf g}\odot[{\bf 3}^3]_{\bf 7}.{\bf 1}\]
\[+2^5.3^4.5.7.[{\bf g}\odot{\bf g}]_{\bf 9}\odot[{\bf 3}^4]_{\bf 9}-\frac{2^8.3^3}{7}\,([{\bf g}\odot{\bf g}]_{\bf 5}\odot{\bf 3}^2)[{\bf 3}^2]_{\bf 1}\]
\[+2^4.3^3[{\bf g}\odot{\bf g}]_{\bf 1}[{\bf 3}^2]_{\bf 1}^2=0.\]

\vspace{.1in}

For economy of notation we have written ${\bf 3}^n$
for the symmetric $n$-fold tensor product of $\bf 3$ and taken for
granted the projection onto the one dimensional (invariant)
component throughout. We have used the identity ${\bf 1}=-4.{\bf
3}^2_{\bf 1}$ in deriving the above. This could be used again to
write the Kummer as an inhomogeneous sextic in the three variables
${\bf 3}$ alone.

This expression for the Kummer should be compared to that in
\cite{CF}. It has the same structure, it being understood that the
variables are differently defined.

\section{Conclusions}

Moving up to higher genus in any approach to the algebraic Jacobian
would appear to be an insane endeavour although an equivariant
description of the Kummer for higher genus might be informative.

However a few specific points deserve to be pursued.

Is there a reason for the dimensional breakdown of the quadratic
identities? The seventy two identities decompose as

\[{\bf
1}^{\oplus 3}\oplus{\bf 2}^{\oplus 4}\oplus{\bf 3}^{\oplus 4}\oplus
{\bf 4}^{\oplus 4}\oplus{\bf 5}^{\oplus 4}\oplus{\bf 6}\oplus{\bf
7}.\]

What is the reason for the sequence of multiplicities
$(3,4,4,4,4,1,1,0,\dots)$?

A related question is why the tensor products of the coordinate
modules have the pole structures they do, related as they are to the
dimensions of the irreducible components in a not quite linear way.
This kind of structure is also apparent in analytic treatments of
the Jacobian via generalised $\wp$-functions \cite{EA} where the
Hirota derivative plays an important r\^ole.

This in turn leads to the question of the relation of the algebraic
approach to the equivariant Kleinian approach of \cite{Ath1,Ath2}.
It is implicit in the definition of the equivariant $\wp$-functions.
Putting $\wp=(\wp_{22},\wp_{12},\wp_{11}))$ this is, in the current
notation, \[[{\bf 3}\odot\wp]_{\bf 1}={\bf 1}.\]

\bibliographystyle{amsplain}

\end{document}